\documentclass[12pt]{article}
\usepackage[expansion=false]{microtype}
\usepackage[reqno]{amsmath}
\usepackage{fixmath}
\usepackage{amssymb, amsthm, enumitem}
\usepackage{mathtools}
\usepackage{graphicx}

\usepackage{ifthen,amssymb,amsfonts,amsmath,amsthm,graphicx, color, tikz, pgf}

\usepackage{hyperref}% wants to be last
\DeclarePairedDelimiter{\norm}{\lVert}{\rVert}

\newtheoremstyle{plainsl}%
    {\topsep}
    {\topsep}
    {\slshape} % only non-default setting
    {}
    {\normalfont\bfseries}
    {.}
    { }
    {}

\swapnumbers

{\theoremstyle{plainsl}
\newtheorem{thm}{Theorem}[section]
\newtheorem{lem}[thm]{Lemma}
\newtheorem{cor}[thm]{Corollary}
\newtheorem{prop}[thm]{Proposition}}
{\theoremstyle{remark}
\newtheorem{ex}[thm]{Example}
}

\renewcommand\proof{\noindent\textsl{Proof. }}
\newcommand\sqr[2]{{\vbox{\hrule height.#2pt
   \hbox{\vrule width.#2pt height#1pt \kern#1pt
        \vrule width.#2pt}\hrule height.#2pt}}}
\renewcommand\qed{%
    \ifmmode\eqno\sqr53
    \else\nolinebreak\ \hfill\sqr53\medbreak\fi}

\numberwithin{equation}{section}

\newcommand\pmat[1]{\begin{pmatrix} #1 \end{pmatrix}}

\newcommand{\abs}[1]{\left | #1 \right |}
\newcommand{\lp}{\! \left (}
\newcommand{\rp}{\right )} \newcommand{\lb}{\left [}
  \newcommand{\rb}{\right ]} \newcommand{\lsb}{\left \{ }
\newcommand{\rsb}{\right \} }

\newcommand{\ver}[1]{\mathbf{e}_{ #1}}
\newcommand\ip[2]{\langle#1,#2\rangle}

\newcommand{\verp}[1]{\mathbf{E}_{ #1}}

\newcommand{\vtxa}{u}
\newcommand{\vtxb}{v}
\newcommand{\vtxc}{w}

\title{Polynomial Characterizations of Distance-Biregular Graphs}

 \author{Sabrina Lato}

\begin{document}
\maketitle

\begin{abstract}Fiol, Garriga, and Yebra~\cite{garriga} introduced the notion of pseudo-distance-regular vertices, which they used to come up with a new characterization of distance-regular graphs. Building on that work, Fiol and Garriga~\cite{garriga} developed the spectral excess theorem for distance-regular graphs. We extend both these characterizations to distance-biregular graphs and show how these characterizations can be used to study bipartite graphs with distance-regular halved graphs and graphs with the spectrum of a distance-biregular graph.
\end{abstract}

\section{Introduction}
Distance-regular graphs are a class of graphs with nice algebraic and structural properties. One definition is that a graph \( G \) of diameter \( d \) with adjacency matrix \( A \) is distance-regular if \( G \) has \( d+1 \) distinct eigenvalues and there exist a sequence of polynomials \( p_0, \ldots, p_d \) where \( p_i \) has degree \( i \) and \( p_i \lp A \rp \) is the \( i \)-th distance adjacency matrix. Distance-regular graphs are well-studied, and more information can be found in the book by Brouwer, Cohen, and Neumaier~\cite{bcn} or the more recent survey by Van Dam, Koolen, and Tanaka~\cite{vanDistance}. In this paper, we will be primarily concerned with two characterizations of distance-regular graphs proven by Fiol, Garriga, and Yebra~\cite{yebra} and Fiol and Garriga~\cite{garriga}.

There is a natural extension of distance-regular graphs to semiregular graphs which gives the class of distance-biregular graphs. If \( G = \lp B, C \rp \) is a bipartite graph with bipartite sets \( B, C,\) adjacency matrix \( A \), and diameter \( d \), then \( G \) is distance-biregular if there exist two sequences of polynomials \( p_0^{B}, \ldots, p_d^{B} \) and \( p_0^{C}, \ldots, p_d^{C} \) such that, for \( S = B, C\) the polynomial \( p_i^{S} \) has degree \( i \) and the restriction of \( p_i \lp A \rp \) to the rows indexed by vertices in \( S \) is the same as the \( i \)-th distance adjacency matrix restricted to the vertices in \( S \). Distance-biregular graphs have been studied by Delorme~\cite{delorme} and Godsil and Shawe-Taylor~\cite{distanceRegularised}, and contain a number of interesting classes of graphs including bipartite distance-regular graphs, generalized polygons, the incidence graphs of partial geometries, the incidence graphs of certain kinds of quasisymmetric designs, points and hyperplanes of affine spaces, and bipartite analogues to the Johnson and Grassmann graphs. 

Fiol, Garriga, and Yebra~\cite{yebra} proved that a graph of diameter \( d \) with adjacency matrix \(  A\) is distance-regular if and only if \( G \) has \( d+1 \) eigenvalues and there exists a polynomial \( p \) of degree \( d \) such that \( p \lp A \rp \) is the \( d \)-th distance adjacency matrix. Fiol~\cite{fiolAlgebraic} included that in a list of characterizations of bipartite distance-regular graphs, and asked for extensions to distance-biregular graphs. In this paper, we provide such an extension by showing that a bipartite graph \( G = \lp B, C \rp \) of diameter \( d \) is distance-biregular if and only if \( G \) has \( d+1 \) distinct eigenvalues and there exist polynomials \( p^{B}, p^{C} \) of degree \( d \) such that \( p^{B} \lp A \rp \) and \( p^C \lp A \rp \) restricted to the rows indexed by vertices in, respectively, \( B \) and \( C \), is the same as the \( d \)-th distance matrix restricted to the rows indexed by vertices in \( B \) or \( C \). 

Building on that work, Fiol and Garriga~\cite{garriga} proved that a regular graph \( G \) with diameter \( d \) and largest eigenvalue \( \lambda \) is distance-regular if and only if there exists a polynomial \( p \) of degree \( d \) such that for every vertex \( \vtxa, \) 
\[ \norm{p}_G^2 = p \lp \lambda \rp = \abs{\lsb \vtxb : d \lp \vtxa, \vtxb \rp = d \rsb}, \]
where \( \norm{}_G \) is a specific norm determined by the spectrum. This theorem is sometimes known as the spectral excess theorem for distance-regular graphs.

We extend their results to bipartite semiregular graphs by defining norms \( \norm{}_B, \norm{}_{C} \) which are determined by the spectra of the graph. Then a semiregular bipartite graph \( G= \lp B, C \rp \) with diameter \( d \) and largest eigenvalue \( \lambda \) is distance-biregular if and only if \( G \) has \( d+1 \) distinct eigenvalues and there exist polynomials \( p^{B}, p^{C} \) of degree \( d \) such that
\[ \norm{p^{B}}_B^2 = p^{B} \lp \lambda \rp = \abs{\lsb \vtxa : d \lp \vtxb, \vtxa \rp = d \rsb} \]
for all \( \vtxb \in B \) and
\[ \norm{p^{C}}_C^2 = p^{C} \lp \lambda \rp = \abs{\lsb \vtxa : d \lp \vtxc, \vtxa \rp = d \rsb} \]
for all \( \vtxc \in C. \) This gives a spectral excess theorem for distance-biregular graphs. Another spectral excess theorem was proven by Fiol~\cite{spectralExcess}, though the conditions there are more involved.

The main idea of the characterizations in~\cite{yebra} and~\cite{garriga} is the use of pseudo-distance-regular vertices. This notion was first introduced by Fiol, Garriga, and Yebra~\cite{yebra}, and it is a local property that vertices in a graph can have. When every vertex in a graph is pseudo-distance-regular, Fiol~\cite{pseudoGlobal} proved that the graph is distance-regular or distance-biregular. Thus the general approach of these characterizations is to come up with conditions on the graph that force every vertex to be pseudo-distance-regular and thereby show that we have a characterization for distance-regular or distance-biregular graphs. Pseudo-distance-regular vertices were studied further by Fiol~\cite{pseudo}.

For a bipartite graph, the halved graphs are the connected components of the distance-two graph. If a graph \( G \) is distance-biregular, then the spectrum of the halved graphs will be determined by the spectrum of \( G \), and the halved graphs will be distance-biregular. We look at when the converse is true. We provide an example from Delorme~\cite{delorme} and a new example of regular bipartite graphs that are not distance-biregular but have distance-regular halved graphs with the spectrum determined from the spectrum of the bipartite graph. Using the characterizations from this paper, we are then able to prove that the only such graphs will be regular. The example of Delorme~\cite{delorme} is a regular bipartite graph of diameter \( d \) that is not distance-regular, even though the \( \lp  d-1 \rp \)-th distance matrix can be written as a polynomial in the adjacency matrix of degree \( d-1 \). This resolves a question posed by Fiol~\cite{fiolAlgebraic}.

One use of the spectral excess theorem for distance-regular graphs that Fiol and Garriga~\cite{garriga} proved is in characterizing when a graph with the spectrum of a distance-regular graph will be distance-regular. This is a problem that has been studied by Haemers~\cite{haemersSpectrumDistance}, Van Dam and Haemers~\cite{vanSpectral}~\cite{whichGraphs}, and Abiad, Van Dam, and Fiol~\cite{abiad}. We conclude this paper by using the spectral excess theorem for distance-biregular graphs presented here to study when a graph which is cospectral to a distance-biregular graph is itself distance-biregular. 

\section{Preliminaries}
Let \( G \) be a connected graph of diameter \( d \) with distance adjacency matrices \( I = A_0, A = A_1, A_2, \ldots, A_d \). Let \( \vtxa \) be a vertex of \( G \).

As described in sources such as Section 8.12 such as Godsil and Royle~\cite{yellow}, the adjacency matrix \( A \) admits a spectral decomposition
\[ A = \sum_{\theta} \theta E_{\theta}, \]
for eigenvalues \( \theta \) and corresponding pairwise orthogonal spectral idempotents \( E_{\theta}. \)
Letting \( \ver{\vtxa} \) be the characteristic vector for \( \vtxa, \) we define the \textit{eigenvalue support} of \( \vtxa \) as the set
\[ \Phi_{\vtxa} = \lsb \theta : E_{\theta} \ver{\vtxa} \neq \mathbf{0} \rsb. \]

The same way that the number of distinct eigenvalues is greater than the diameter of a graph, the size of the eigenvalue support of a vertex is greater than its eccentricity, a result that can be found, for example, in Section 5.2 of Coutinho and Godsil~\cite{Godsil2017a}.

\begin{lem}\label{eccentricSupport}Let \( \vtxa \) be a vertex in graph \( G \) with eccentricity \( e \). Then
  \[ \abs{\Phi_{\vtxa}} \geq e+1. \]
\end{lem}

\proof Let \( U \) be the cyclic \( A \)-module generated by \( \ver{\vtxa} \). The eigenvectors corresponding to distinct eigenvalues in the eigenvalue support form a basis for \( U \). Let \( 0 \leq r \leq e, \) and let \( \vtxb \) be at distance \( r \) from \( \vtxa \). Then
\[ \ver{\vtxb} A^r \ver{\vtxa} \neq 0, \]
but for all \( s < r, \) we have
\[ \ver{\vtxb} A^s \ver{\vtxa} = 0. \]
Thus \( A^r \ver{\vtxa} \) cannot be expressed as a linear combination of \( A^0 \ver{\vtxa}, \ldots, A^{r-1} \ver{\vtxa}, \) and so \( \ver{\vtxa}, A \ver{\vtxa}, \ldots, A^e \ver{\vtxa} \) are linearly independent, and so \( U \) has dimension at least \( e+1 \).\qed

When the bound is tight, that is, when \( \abs{\Phi_{\vtxa}}=e+1, \) we say that vertex \( \vtxa \) is \textit{spectrally extremal}.

We can define a \( \vtxa \)-local inner product on polynomials associated to by
\[ \ip{p}{q}_{\vtxa} = \ver{\vtxa}^T p \lp A \rp q \lp A \rp \ver{\vtxa}. \]
Using the spectral decomposition, this is equivalent to
\[ \ip{p}{q}_{\vtxa} = \sum_{\theta \in \Phi_{\vtxa}} \ver{\vtxa}^T E_{\theta} \ver{\vtxa} p \lp \theta \rp q \lp \theta \rp. \]

A sequence of polynomials \( \lp p_i \rp_{i \geq 0} \) is a sequence of \textit{orthogonal polynomials} with respect to this inner product if the polynomials are pairwise orthogonal and \( p_i \) has degree \( i \) for \( i \geq 0. \) An easy consequence, which can be found in texts on orthogonal polynomials such as Nikiforov, Suslov, and Uvarov~\cite{classicalOrthogonal}, says that orthogonal polynomials satisfy a particular three-term recurrence. That is, for all \( i \geq 1, \) there exist real numbers \( a_i, b_i, c_i \) such that
\[ x p_i \lp x \rp = c_{i+1} p_{i+1} \lp x \rp + a_i p_{i} \lp x \rp + b_{i-1} p_{i-1} \lp x \rp. \]

A vertex \( \vtxa \) is \textit{locally distance-regular} if, for any vertex \( \vtxb \) at distance \( i \) from \( \vtxa, \) the numbers of vertices adjacent to \( \vtxb \) and at distances \( i-1, i \) and \( i+1 \) are independent of the choice of vertex \( \vtxb. \) We can use orthogonal polynomials to obtain an equivalent definition of locally distance-regular.

\begin{lem}Let \( G \) be a connected graph with a vertex \( \vtxa \) of eccentricity \( e \). Then \( \vtxa \) is locally distance-regular if and only if \( \abs{\Phi_{\vtxa}} = e+1 \) and, for \( 0 \leq i \leq e, \) there exist a sequence of orthogonal polynomials \( p_0, \ldots, p_e \) satisfying
  \[ p_i \lp A \rp \ver{\vtxa} = A_i \ver{\vtxa}. \]
\end{lem}

\proof Let
\[ p_{e+1} = \prod_{\theta_r \in \Phi_{\vtxa}} \lp x- \theta_r \rp, \]
which we append to our sequence of orthogonal polynomials. Then using the three-term recurrence, the existence of such a sequence of polynomials is equivalent to the property that, for all \( 0 \leq i \leq e, \) we have
\[ A A_i \ver{\vtxa} = b_{i-1} A_{i-1} \ver{\vtxa} + a_i A_i \ver{\vtxa} + c_{i+1} A_{i+1} \ver{\vtxa}. \]
Taking the \( \lp \vtxa, \vtxb \rp \)-th entry, this is equivalent to the property that for any \( \vtxb \) and any \( 0 \leq i \leq e, \) we have
\[ \abs{\lsb \vtxc \sim \vtxb : d \lp \vtxa, \vtxc \rp = i \rsb} =  \begin{cases}b_{i-1} & d \lp \vtxa, \vtxb \rp = i-1 \\ a_i & d \lp \vtxa, \vtxb \rp = i \\ c_{i+1} & d \lp \vtxa, \vtxb \rp = i+1 \end{cases}, \]
which is equivalent to \( \vtxa \) being locally distance-regular.\qed

We can extend the local inner product, polynomials, and idea of locally distance-regular vertices to apply to sets of vertices. Let \( S \subseteq V(G) \) be a set of vertices. We can define the eigenvalue support of \( S \) by
\[ \Phi_S = \abs{\bigcup_{s \in S} \Phi_{s}} \]
and the \( S \)-local inner product by
\[ \ip{f}{g}_S = \frac{1}{\abs{S}} \sum_{\vtxa \in S} \ver{\vtxa}^T f \lp A \rp g \lp A \rp \ver{\vtxa}. \]

Let \( \verp{S} \) be the \( \abs{V \lp G \rp} \times \abs{S} \) matrix whose columns are the characteristic vectors for vertices in \( S \). The set \( S \) is locally distance-regular if \( \abs{\Phi_S} -1 \) is the eccentricity of every vertex in \( S \) and there exists a sequence of orthogonal polynomials \( p_0^{S}, \ldots, p^{S}_{\abs{\Phi_S}-1} \) such that, for all \( 0 \leq i \leq \abs{\Phi_S}-1, \) we have
\[ p_i^{S} \lp A \rp \verp{S} = A_i \verp{S}. \]

Godsil and Shawe-Taylor~\cite{distanceRegularised} proved that a graph which is locally distance-regular at every vertex is either distance-regular or distance-biregular. Their results can be interpreted as saying that if there is a partition of the vertices into locally-distance-regular sets, the partition is either the trivial partition or defines a bipartition on the graph.

In the language of polynomials, then, a graph \( G \) of diameter \( d \) is distance-regular if has \( d+1 \) distinct eigenvalues and there exists a sequence of polynomials \( p_0, \ldots, p_d \) such that that for all \( 0 \leq i \leq d, \) the polynomial \( p_i \) has degree \( i \) and satisfies \( p_i \lp A \rp = A_i. \)

Let \( G = \lp B, C \rp \) be a bipartite graph with diameter \( d \). Assume without loss of generality that vertices in \( B \) have maximum eccentricity \( d \). Then \( G \) is distance-biregular if and only if it has \( d+1 \) distinct eigenvalues and there exist sequences of polynomials \( p_0^{B}, p_1^{B}, \ldots, p_d^{B} \) and \( p_0^{C}, p_1^{C}, \ldots, p_d^{C} \) such that \( p_i^{B}, p_i^{C} \) have degree \( i \) and satisfy
\[ p_i^{B} \lp A \rp \verp{B} = A_i \verp{B} \]
and
\[ p_i^{C} \lp A \rp \verp{C} = A_i \verp{C}. \]

When discussing the inner product local to sets, we will generally be talking about the set of all vertices, or the sets of a bipartite graph.

Note that
\[ \ip{p}{q}_{G} = \frac{1}{\abs{V}} \sum_{\theta \in \Phi_{\vtxa}} \sum_{\vtxa \in V} \ver{\vtxa}^T E_{\theta} \ver{\vtxa} p \lp \theta \rp q \lp \theta \rp = \sum_{\theta} \frac{\mathrm{tr} \lp E_{\theta} \rp}{\abs{V}} p \lp \theta_r \rp q \lp \theta_r \rp, \]
and since tr\( \lp E_{\theta} \rp \) is the multiplicity of \( \theta, \) the global inner product is determined by the spectrum.

If \( G = \lp B, C \rp \) is a bipartite graph, then we can write the adjacency matrix as a block matrix of the form
\[ A = \pmat{0 & N \\ N^T & 0}, \]
where \( N \) is the biadjacency matrix from vertices in \( B \) to vertices in \( C \). Then for any nonzero eigenvalue \( \theta, \) we can write the spectral idempotent as a block matrix of the form
\[ \pmat{B_{\theta} & D_{\theta} \\ D_{\theta}^T & C_{\theta}}, \]
and observe
\[ \pmat{B_{\theta} & -D_{\theta} \\ -D_{\theta}^T & C_{\theta}} \]
is the spectral idempotent for \( - \theta. \) This decomposition has a couple of important consequences.

If \( \theta \) is any nonzero eigenvalue of \( G, \) then \( 2 B_{\theta} \) is an idempotent for \( NN^T \) with eigenvalue \( \theta^2 \). Similarly, \( 2 C_{\theta} \) is an idempotent for \( N^TN \). Since \( NN^T \) and \( N^TN \) share nonzero with multiplicities, we have that
\[ \mathrm{tr}\lp B_{\theta} \rp = \mathrm{tr} \lp C_{\theta} \rp = \frac{1}{2} \mathrm{tr} \lp E_{\theta} \rp. \]
Suppose without loss of generality that \( \abs{B} \geq \abs{C} \). Then the multiplicity of \( 0 \) is \( \abs{B}-\abs{C}, \) plus twice the multiplicity of 0 as an eigenvalue for \( N^TN \). Therefore
\[ \mathrm{tr} \lp B_{\theta} \rp - \abs{B}+\abs{C} = \mathrm{tr} \lp C_{\theta} \rp = \frac{1}{2} \lp \mathrm{tr} \lp E_{\theta} \rp -\abs{B}+\abs{C} \rp. \]
If we know that \( G \) is a \( \lp k, \ell \rp \)-semiregular graph, then we can determine \( \abs{B} \) and \( \abs{C} \) from the spectrum. Therefore, for bipartite semiregular graphs, the inner products \( \ip{}{}_B \) and \( \ip{}{}_C \) are also determined by the spectrum.

\begin{lem}\label{sumZeroOff}If \( E_{\theta} \) and \( E_{\tau} \) are two distinct spectral idempotents such that \( E_{\theta}+ E_{\tau} \) is zero on the off-diagonal blocks, then \( E_{\theta} = E_{-\tau}. \)\end{lem}

  \proof We have that \( D_{\theta} + D_{\tau} = 0, \) so for the ease of notation let let \( D_{\theta} = D \) and \( D_{\tau} = -D. \) We may assume without loss of generality that \( \theta \neq 0. \) Then \( \theta \neq -\theta, \) and therefore \( E_{\theta} \) is orthogonal to \( E_{-\theta}. \) Suppose that \( E_{-\theta} \neq E_{\tau}. \) Then \( E_{-\theta} - E_{\tau} \) is a nonzero matrix that is zero on the off-diagonal blocks. It is also orthogonal to \( E_{\theta}, \) so 
\[ B_{\theta} \lp B_{\theta} - B_{\tau} \rp = \mathbf{0}. \]
Since \( E_{-\theta} \) and \( E_{\tau} \) are distinct, they are also orthogonal, and so
\[ \mathbf{0} = B_{\theta} B_{\tau} + DD^T, \]
and since \( E_{\theta} \) is idempotent we have
\[ B_{\theta} = B_{\theta}^2 + DD^T = B_{\theta}^2 - B_{\theta} B_{\tau} = \mathbf{0}. \]
An analogous argument gives us that \( C_{\theta} = \mathbf{0}, \) which is a contradiction, and therefore \( E_{\tau} = E_{-\theta}. \)\qed

\begin{lem}\label{biRegularOddDiameter}Let \( G = \lp B, C \rp \) be a semiregular bipartite graph with diameter \( d \) and \( d+1 \) eigenvalues. If \( d \) is odd, then \( G \) is regular.
\end{lem}

\proof Since \( d \) is odd, \( d+1 \) is even, and therefore zero cannot be an eigenvalue of \( G \). On the other hand, zero is an eigenvalue of \( G \) with multiplicity \( \abs{\abs{B}-\abs{C}} \), so \( \abs{B}= \abs{C} \), and since \( G \) is \( \lp k, \ell \rp \)-semiregular, we can therefore conclude that \( k = \ell, \) and so \( G \) is regular.\qed

Let \( G \) be a connected graph with largest eigenvalue \( \lambda \). By Perron-Frobenius, there exists an eigenvector of norm one for \( \lambda \) that has all positive entries. We will refer to this eigenvector as the \textit{Perron eigenvector}, and denote it \( \mathbf{v}. \) The entry of \( \mathbf{v} \) indexed by \( \vtxa \) is \( v_{\vtxa}. \)

If \( G \) is a regular graph, then \( \mathbf{v} = \frac{1}{\sqrt{\abs{V}}} \mathbf{1}. \)

If \( G \) is a \( \lp k, \ell \rp \)-semiregular graph, then \( \lambda = \sqrt{k \ell} \) and
\[ \mathbf{v} = \frac{1}{k \abs{B} + \ell \abs{C}} \pmat{\sqrt{k} \\ \vdots \\ \sqrt{k} \\ \sqrt{\ell} \\ \vdots \\ \sqrt{\ell}} = \frac{1}{2 \abs{B}} \pmat{1 \\ \vdots \\ 1 \\ \sqrt{\frac{\ell}{k}} \\ \vdots \\ \sqrt{\frac{\ell}{k}}} = \frac{1}{2 \abs{C}} \pmat{\sqrt{\frac{k}{\ell}} \\ \vdots \\ \sqrt{\frac{k}{\ell}} \\ 1 \\ \vdots \\ 1}. \]

\section{Pseudo-Distance-Regular Vertices}\label{pseudo}

We say that a vertex \( \vtxa \) of eccentricity \( e \) is \textit{pseudo-distance-regular} if \( \abs{\Phi_{\vtxa}} = e+1 \) and, for all \( 0 \leq i \leq e, \) there exists a polynomial \( p_i \lp x \rp \) of degree \( i \) such that \( p_i \lp A \rp \ver{\vtxa} \) is nonzero precisely on the vertices at distance \( i \) from \( \vtxa \). Note that such a sequence \( p_0, \ldots, p_e \) is clearly orthogonal with respect to the \( \vtxa \)-local inner product.

Since \( \vtxa \) is spectrally extremal, we have two basis for the cyclic \( A \)-module generated by \( \vtxa \). The first is given by eigenvectors corresponding to eigenvalues in the eigenvalue support of \( \vtxa, \) the second is given by \( p_0 \lp A \rp \ver{\vtxa}, \ldots, p_e \lp A \rp \ver{\vtxa}. \)

Using the second basis, we can write the Perron eigenvector \( \mathbf{v} \) as a linear combination of \( p_i \lp A \rp \ver{\vtxa} \) for \( 0 \leq i \leq e. \) Since the entries of \( \mathbf{v} \) are positive, and the vectors \( p_i \lp A \rp \ver{\vtxa} \) are nonzero on disjoint entries, the coefficients of the linear combination must be nonzero. In particular, the nonzero entries of \( p_i \lp A \rp \) must be some scalar multiple of the entries of \( \mathbf{v} \) indexed by vertices at distance \( i \) from \( \vtxa \). 

Pseudo-distance-regular vertices were introduced by Fiol, Garriga, and Yebra~\cite{yebra}, using a different definition. However, the definition used here is, by Theorem 3.10 of~\cite{yebra}, equivalent.

The existence of the last pseudo-distance polynomial is sufficient to define pseudo-distance-regularity, as was shown in Theorem 6.3 of~\cite{yebra}. We will reprove this result here, with a different set-up.

\begin{thm}\label{pseudoLast}Let \( G \) be a connected graph with vertex \( \vtxa \) of eccentricity \( e \). Then \( G \) is pseudo-distance-regular at \( \vtxa \) if and only if \( \abs{\Phi_{\vtxa}} = e+1 \) and there exists a polynomial \( p \) of degree \( e \) such that \( p \lp A \rp \ver{\vtxa} \) is nonzero precisely on the vertices at distance \( e \) from \( \vtxa \).
\end{thm}

\proof By definition, if \( G \) is pseudo-distance-regular at \( \vtxa \) then \( p_e \lp A \rp \ver{\vtxa} \) has the desired property. Thus we may let \( \vtxa \) be a spectrally extremal vertex with eccentricity \( e \) and suppose there exists a polynomial \( p \) of degree \( e \) such that \( p \lp A \rp \ver{\vtxa} \) is nonzero precisely on the vertices at distance \( e \) from \( \vtxa \).

Let \( p_0, \ldots, p_{e-1} \) be a sequence of orthogonal polynomials. If  \( 0 \leq i < e \) then for any vertex \( \vtxb \) at distance \( e \) from \( \vtxa, \) we have \( \ver{\vtxb}^T A^i \ver{\vtxa} = 0, \) and therefore
\[ \ip{p}{p_i} = \ver{\vtxa}^T p \lp A \rp p_i \lp A \rp \ver{\vtxa} = 0, \]
Since \( p \) is orthogonal to \( p_0, \ldots, p_{e-1}, \) we may append it to our sequence of polynomials as \( p_e. \) We also append 
\[ p_{e+1} \lp x \rp = \prod_{\theta_r \in \Phi_{\vtxa}} \lp x-\theta_r \rp. \]
to our sequence.

It is true by our definitions that the only nonzero entries of \( p_{e+1} \lp A \rp \ver{\vtxa} \) and \( p_e \lp A \rp \ver{\vtxa} \) are indexed by vertices at distance, respectively, \( e+1 \) and \( e \) from \( \vtxa \). Proceeding inductively, we may assume that the only nonzero entries of \( p_i \lp A \rp \ver{\vtxa}, p_{i+1} \lp A \rp \ver{\vtxa}, \ldots, p_{e+1} \lp A \rp \ver{\vtxa} \) are indexed by vertices at the right distance from \( \vtxa \). Since \( p_0, p_1, \ldots, p_{e+1} \) is a sequence of orthogonal polynomials, we know that there exist values \( a_i, b_{i-1}, c_{i+1} \) such that
\[ x p_i \lp x \rp = b_{i-1} p_{i-1} \lp x \rp + a_i p_i \lp x \rp + c_{i+1} p_{i+1} \lp x \rp. \]
In particular, this gives us
\[ A p_i \lp A \rp \ver{\vtxa} = b_{i-1} p_{i-1} \lp A \rp \ver{\vtxa} + a_i p_i \lp A \rp \ver{\vtxa} + c_{i+1} p_{i+1} \lp A \rp \ver{\vtxa}. \]
On the other hand,
\[ \ver{\vtxb}^T A p_i \lp A \rp \ver{\vtxa} = \sum_{\vtxc \sim \vtxb} \ver{\vtxc} p_i \lp A \rp \ver{\vtxa}.\]
By the inductive hypothesis, the nonzero entries of \( p_{i+1} \lp A \rp \ver{\vtxa} \) and \( p_{i} \lp A \rp \ver{\vtxa}, \) are indexed by vertices at distance, respectively \( i+1 \) and \( i \) from \( \vtxa \). Thus,we must have that the nonzero entries of \( p_{i-1} \lp A \rp \ver{\vtxa} \) can only be at distance \( i-1 \) from \( \vtxa \). Since \( p_{i-1} \) is a polynomial of degree \( i-1, \) we conclude \( p_{i-1} \lp A \rp \ver{\vtxa} \) is nonzero precisely on the vertices at distance \( i- 1 \) from \( \vtxa. \)\qed

As we did for locally distance-regular vertices, we can extend this notion to sets. The previous proof extends naturally in this context.

\begin{cor}\label{setPseudo}Let \( G \) be a connected graph, and let \( S \) be a set of vertices which all have eccentricity \( e \). Then \( G \) is pseudo-distance-regular at \( S \) if and only if \( \abs{\Phi_S} = e+1 \) and there exists a polynomial \( p_e \) of degree \( e \) such that \( p_e \lp A \rp \verp{S} \) is nonzero precisely on the vertices at distance \( e \) from \( \vtxa \).
\end{cor}

Godsil and Shawe-Taylor~\cite{distanceRegularised} proved that a graph where every vertex was locally distance-regular was either distance-regular or distance-biregular. Fiol, Garriga, and Yebra~\cite{yebra} showed that if every vertex was pseudo-distance-regular, then every vertex was locally distance-regular and therefore the graph was distance-regular or distance-biregular. Subsequently, Fiol~\cite{pseudoGlobal} gave an independent proof of this.

\begin{thm}\label{pseudoGlobal}[Fiol~\cite{pseudoGlobal}]Let \( G \) be a connected graph that is pseudo-distance-regular at every vertex. Then \( G \) is either distance-regular or distance-biregular.
\end{thm}

It is not difficult to determine when a graph will be distance-regular and when it will be distance-biregular. Indeed, any regular graph that is pseudo-distance-regular at every vertex is distance-regular, and any bipartite graph that is pseudo-distance-regular at every vertex is distance-biregular. Any bipartite regular graph which is pseudo-distance-regular at every vertex is both distance-regular and distance-biregular.

In the parlance of pseudo-distance-regular sets, a graph with vertex set \( V \) is distance-regular when the set \( V \) is pseudo-distance-regular. Similarly, a bipartite graph \( G = \lp B, C \rp\) is distance-biregular when \( B \) and \( C \) are both pseudo-distance-regular sets.

\section{Diametral Polynomial Characterization}

One of the major results in Fiol, Garriga, and Yebra~\cite{yebra} was a new characterization of when a graph is distance-regular. They showed that a graph \( G \) with diameter \( d \) is distance-regular if and only if \( G \) has \( d+1 \) distinct eigenvalues, every vertex has eccentricity \( d \), and there exists a polynomial \( p \) of degree \( d \) such that \( p \lp A \rp = A_d. \) In fact, we prove that the condition that every vertex has eccentricity \( d \) is not necessary.

\begin{thm}\label{weakDistance}Let \( G \) be a connected graph with diameter \( d \). Then \( G \) is distance-regular if and only if \( G \) has \( d+1 \) distinct eigenvalues and there exists a polynomial \( p \) of degree \( d \) such that \( p \lp A \rp = A_d. \)
\end{thm}

\proof A distance-regular graph clearly has the desired properties, so assume that \( G \) has \( d+1 \) distinct eigenvalues and a polynomial \( p \) of degree \( d \) exists such that \( p \lp A \rp = A_d. \)

We first prove that every vertex in \( G \) has eccentricity \( d \). Suppose otherwise, and let \( \vtxa \) be a vertex of eccentricity \( d-1 \). Then we have
\[ \mathbf{0} = A_d \ver{\vtxa} = p \lp A \rp \ver{\vtxa} = \sum_{\theta \in \Phi{\vtxa}} p \lp \theta \rp E_{\theta} \ver{\vtxa}, \]
and since \( E_{\theta} \ver{\vtxa} \) are eigenvectors for distinct eigenvalues, they must be linearly independent. Thus every eigenvalue in \( \Phi_{\vtxa} \) is a root of \( p. \) By Lemma~\ref{eccentricSupport} we know that \( \abs{\Phi_{\vtxa}} \geq d, \) and since \( p \lp A \rp \) is nonzero, there must be precisely one eigenvalue \( \hat{\theta} \) of \( G \) which is not a root of \( p \). Thus we have
\[ p \lp A \rp = p \lp \hat{\theta} \rp E_{\hat{\theta}}, \]
and for every vertex \( \vtxb, \) this gives us
\[ p \lp \hat{\theta} \rp \ver{\vtxb}^T E_{\hat{\theta}} \ver{\vtxb} = \ver{\vtxb}^T p \lp A \rp \ver{\vtxb} = \ver{\vtxb}^T A_d \ver{\vtxb} = 0. \]
Since \( p \lp \hat{\theta} \rp \neq 0, \) this implies that \( \ver{\vtxb}^T E_{\hat{\theta}} \ver{\vtxb} = 0. \) Further, because the spectral idempotents are positive semidefinite, this means that \( E_{\hat{\theta}} \ver{\vtxb} = \mathbf{0} \) and thus \( \hat{\theta} \notin \Phi_{\vtxb} \) for any vertex \( \vtxb \). However, since \( G \) has diameter \( d, \) there must be a vertex with eccentricity \( d, \) and by Lemma~\ref{eccentricSupport}, this vertex must have all \( d+1 \) eigenvalues in its support. This gives us a contradiction, and so every vertex in \( G \) must have eccentricity \( d \).

Let \( \vtxa \) be an arbitrary vertex in \( G \). Since \( \vtxa \) has eccentricity \( d \) and \( G \) only has \( d+1 \) distinct eigenvalues, we see that \( \vtxa \) is spectrally extremal. By assumption, we have
\[ p \lp A \rp \ver{\vtxa} = A_d \ver{\vtxa}, \]
and so by Theorem~\ref{pseudoLast}, we see \( \vtxa \) is pseudo-distance-regular. Since our choice of \( \vtxa \) was arbitrary, Corollary~\ref{setPseudo} tells us that the set of all vertices \( V \) is pseudo-distance-regular. Therefore, by Theorem~\ref{pseudoGlobal}, we know that \(G \) is distance-regular.\qed

We can also extend these ideas to get a distance-biregular analogue of Theorem~\ref{weakDistance}. In this case, we have one polynomial for each set of the partition, and we show that, when restricting to the sets of the partition, these polynomials give us the distance matrices. This provides a new characterization of when a bipartite graph is distance-biregular.

\begin{thm}\label{pseudoBiregular}Let \( G = \lp B, C \rp \) be a connected bipartite graph with diameter \( d \). Then \( G \) is distance-biregular if and only if \( G \) has \( d+1 \) distinct eigenvalues and there exist polynomials \( p^{B} \) and \( p^{C} \) of degree \( d \) such that
  \[ p^{B} \lp A \rp \verp{C} = A_d \verp{B} \]
  and
  \[ p^{C} \lp A \rp \verp{C} = A_d \verp{C}. \]
\end{thm}

\proof These conditions clearly hold if \( G \) is distance-biregular, so we may assume that \( G \) is a graph with \( d+1 \) distinct eigenvalues and polynomials \( p^{B}, p^{C}. \) Assume without loss of generality that \( d \) is the maximum eccentricity of vertices in \( B \), and let \( d' \) be the maximum eccentricity of vertices in \( C \). We begin by showing that every vertex in \( B \) has eccentricity \( d \) and every vertex in \( C \) has eccentricity \( d' \). If \( d' = d, \) then the same argument as in the proof of Theorem~\ref{weakDistance} applies, and so every vertex has eccentricity \( d \). Thus, we can assume \( d' = d-1 \) and \( d \) is even. We wish to show that every vertex in \( B \) has eccentricity \( d \).

Writing \( p^{B} \lp A \rp \) as a block matrix, we see that \( p^{B} \lp A \rp \) is zero on the off-diagonal blocks, and therefore \( p^{B} \) is an even function. 

Suppose there exists a vertex \( \vtxb \in B \) with eccentricity \( d-2. \) Then since
\[  \mathbf{0} = A_d \ver{\vtxb} = p^{B} \lp A \rp \ver{\vtxb}, \]
we again see that every eigenvalue in \( \Phi_{\vtxa} \) is a root of \( p^{B}. \) Since \( p^{B} \lp A \rp \) is nonzero, there must be at most two eigenvalues \( \theta, \tau \) of \( G \) which are not roots of \( p^{B}. \) Then
\[ p^{B} \lp A \rp = p^{B} \lp \theta \rp E_{\theta} + p^{B} \lp \tau \rp E_{\tau}. \]
Since the off-diagonal blocks of \( p^{B} \lp A \rp \) are zero, Lemma~\ref{sumZeroOff} tells us that \( \tau = - \theta. \) 

Let \( \vtxa \in B \) have eccentricity \( d \). Then all \( d+1 \) eigenvalues of \( G \) must be in the support of \( \vtxa \). On the other hand, we have
\[ 0 = \ver{\vtxa}^T p^{B} \lp A \rp \ver{\vtxb} = p^{B} \lp \theta \rp \ver{\vtxa}^T E_{\theta} \ver{\vtxa} + p^{C} \lp -\theta \rp \ver{\vtxa}^T E_{-\theta} \ver{\vtxa}, \]
and since every eigenvalue of \( G \) is in \( \Phi_{\vtxa} \), we must have
\[ \ver{\vtxa}^T E_{\theta} \ver{\vtxa}= \ver{\vtxa}^T E_{-\theta} \ver{\vtxa} \neq 0. \]
Thus \( p^{B} \lp \theta \rp = - p^{B} \lp -\theta \rp, \) which is a contradiction since \( p^{B} \) is an even function. 

Therefore every vertex in \( B \) has eccentricity \( d, \) from which it immediately follows that every vertex in \( C \) has eccentricity \( d'. \)

Since there are only \( d+1 \) distinct eigenvalues, \( \Phi_{B} = d+1, \) and
\[ p^{B} \lp A \rp \verp{B} = A_d \verp{B}, \]
so by Corollary~\ref{setPseudo}, the set \( B \) is pseudo-distance-regular. If \( d' = d, \) an analogous argument shows \( C \) is pseudo-distance-regular. Otherwise, \( d' = d \). Because \( p^{C} \verp{C} = \mathbf{0}, \) we see that \( \abs{\Phi_C} = d = d'+1. \) Further, since we have already shown that the set \( B \) is pseudo-distance-regular, we know that there exists a polynomial \( p_{d-1}^{B} \) of degree \( d-1 \) such that \( p_{d-1}^{B} \lp A \rp \verp{B} \) is nonzero precisely on the vertices at distance \( d-1 \) from vertices in \( B \). But \( p_{d-1}^{B} \) is symmetric, so \( p_{d-1}^{B} \lp A \rp \verp{C} \) must also be nonzero precisely on the vertices at distance \( d-1 \) from vertices in \( C \). Then Corollary~\ref{setPseudo} tells us \( C \) is pseudo-distance-regular. It follows from Theorem~\ref{pseudoGlobal} that \( G \) is distance-biregular.\qed

\section{Distance-Regular Halved Graphs}

Every distance-biregular graph has distance-regular halved graphs. In fact, a stronger statement is true. Let \( G = \lp B, C \rp \) be a \( \lp k, \ell \rp \)-semiregular graph. Let \( H_B \) be the halved graph of vertices in \( B \), and let \( H_C \) be the halved graph of vertices in \( C \). Let \( N \) be the biadjacency matrix from \( B \) to \( C \). If \( G \) is distance-biregular, then \( H_B \) and \( H_C \) are distance-regular and there exist constants \( r,s \) such that
\[ NN^T = r A \lp H_B \rp + k I \]
and
\[ N^TN = s A \lp H_C \rp + \ell I. \]

Spectrally, the latter condition means that the eigenvalues of the halved graphs are determined by the eigenvalues of the bipartite graph and vice-versa. Suppose without loss of generality that \( \abs{B} \geq \abs{C}. \) If \( G \) is distance-biregular with halved graphs \( H_B, H_C, \) then \( G \) is \( \lp k, \ell \rp \)-semiregular graph and there exist constants \( r, s \) such that the eigenvalues of \( G \) are
\[ \lsb \pm \sqrt{r \theta + k} : \theta \in \text{ev} \lp H_1 \rp \rsb = \lsb \pm \sqrt{s \tau + \ell} : \tau \in \text{ev} \lp H_2 \rp \rsb \cup \lsb 0^{(\abs{B}-\abs{C})} \rsb. \]

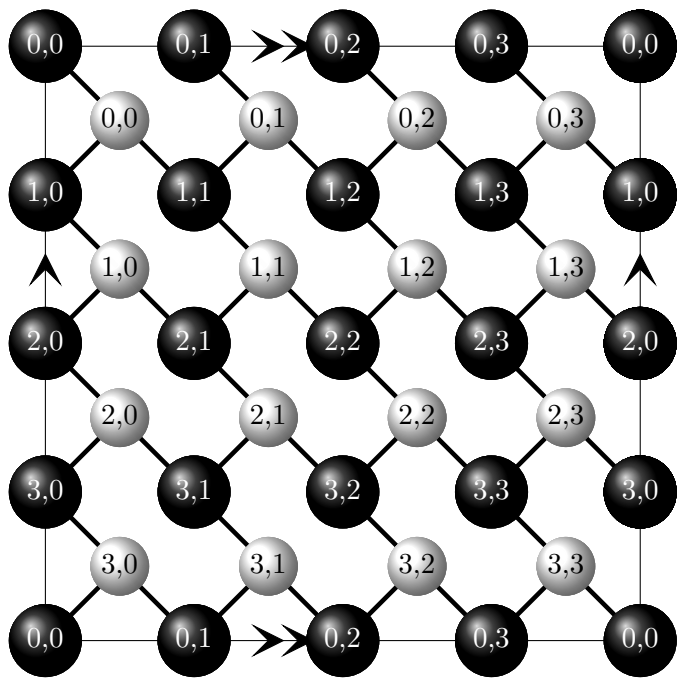
\begin{figure}
  \begin{center}
    \begin{tikzpicture}[x=0.75pt,y=0.75pt,yscale=-.75,xscale=.75, node distance= 1.5cm]
      %uncomment if require: \path (0,420); %set diagram left start at 0, and has height of 420
		
      %Shape: Square [id:dp7863943622357719] 
      \draw   (130,10) -- (530,10) -- (530,410) -- (130,410) -- cycle ;
      \draw   [fill={rgb, 255:red, 0; green, 0; blue, 0 }  ,fill opacity=1 ] (120,170) -- (130,150) -- (140,170) -- (130,160) -- cycle ;
      \draw   [fill={rgb, 255:red, 0; green, 0; blue, 0 }  ,fill opacity=1 ] (520,170) -- (530,150) -- (540,170) -- (530,160) -- cycle ;
      \draw  [fill={rgb, 255:red, 0; green, 0; blue, 0 }  ,fill opacity=1 ] (270,0) -- (290,10) -- (270,20) -- (280,10) -- cycle ;
      \draw  [fill={rgb, 255:red, 0; green, 0; blue, 0 }  ,fill opacity=1 ] (290,0) -- (310,10) -- (290,20) -- (300,10) -- cycle ;
      \draw  [fill={rgb, 255:red, 0; green, 0; blue, 0 }  ,fill opacity=1 ] (270,400) -- (290,410) -- (270,420) -- (280,410) -- cycle ;
      \draw  [fill={rgb, 255:red, 0; green, 0; blue, 0 }  ,fill opacity=1 ] (290,400) -- (310,410) -- (290,420) -- (300,410) -- cycle ;

      \foreach \y in {0,1,2,3}{
        \foreach \x in {0,1,2,3}{
	  \draw [line width=1.75]  (180+100*\y,60+100*\x) -- (130+100*\y,10+100*\x) ;
	  \draw [line width=1.75]  (180+100*\y,60+100*\x) -- (130+100*\y,110+100*\x) ;
	  \draw [line width=1.75]  (180+100*\y,60+100*\x) -- (230+100*\y,110+100*\x) ;          
        }
      }    
      
      \tikzset{VertexStyle/.style = {shape= circle, ball color= white, inner sep= 2pt, outer sep= 0pt, minimum size = 20 pt}}
      
      \foreach \y in {0,1,2,3}{
        \foreach \x in {0,1,2,3}{
        \draw (180+100*\y,60+100*\x) node[VertexStyle, font=\small] (\x-\y) {\x,\y} ;
        }
      }

      \tikzset{VertexStyle/.style = {shape= circle, ball color= black, minimum size = 20 pt}}
      \foreach \y in {0,1,2,3}{
        \foreach \x in {0,1,2,3}{
          \draw (130+100*\y,10+100*\x) node[VertexStyle, font=\small,text=white] (\x+\y) {\x,\y} ;
          \draw (530,10+100*\x) node[VertexStyle, font=\small,text=white] (\x+4) {\x,0} ;
        }
        \draw (130+100*\y,410) node[VertexStyle, font=\small,text=white] (4+\y) {0,\y} ;
      }
      \draw (530,410) node[VertexStyle, font=\small,text=white] (4+4) {0,0} ;     
    \end{tikzpicture}
  \end{center}
  \caption{Example of Delorme~\cite{delorme}} \label{delormeShrikhande}
\end{figure}

The converse does not hold, as was shown by Delorme~\cite{delorme}.

\begin{ex}Consider the graph in Figure~\ref{delormeShrikhande}. This is regular bipartite graph with diameter 5 and spectrum
\[ \lsb 3^{(1)}, \sqrt{5}^{(6)}, 1^{(9)}, -1^{(9)}, -\sqrt{5}^{(6)}, -3^{(1)} \rsb. \]
The halved graphs are the Shrikhande graphs, and the eigenvalues have the necessary relationship. However, there are two walks of length three from the black vertex \( \lp 2, 2 \rp \) to the white vertex \( \lp 2, 1 \rp \), but only one walk of length three from black \( \lp 2, 2 \rp \) to white \( \lp 3, 3 \rp \), so \( A_3 \) cannot be polynomial in \( A, \) and therefore \( G \) is not distance-biregular.
\end{ex}

It is also worth noting that for this graph,
\[ p_4 \lp x \rp := \frac{1}{2} x^4 -4 x^2 + \frac{9}{2} \]
is the fourth-degree polynomial of our orthogonal sequence, and it satisfies \( p_4 \lp A \rp = A_4. \) This resolves Problem 3 in Fiol~\cite{fiolAlgebraic}, since it is a regular bipartite graph with \( p_4 \lp A \rp = A_4 \) that is not distance-regular.

Even though Delorme's example is not distance-biregular, every vertex is still spectrally extremal. However, this is not necessary for a graph with distance-regular halved graphs.

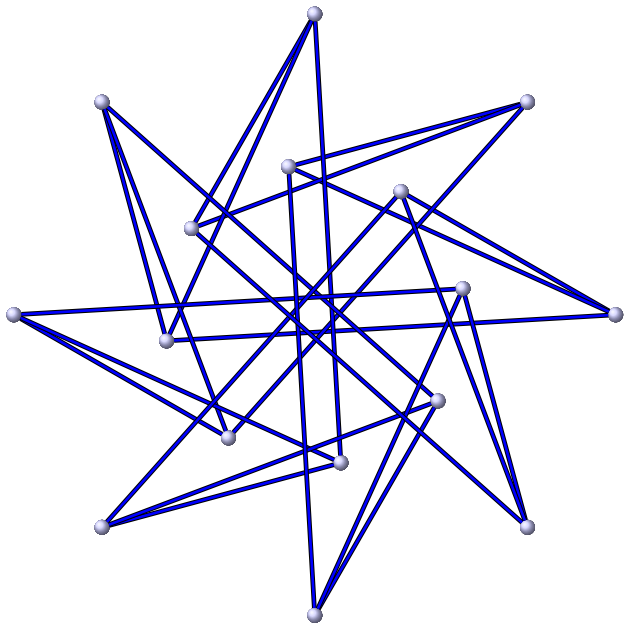
\begin{figure}
  \begin{center}
    \begin{tikzpicture}[node distance = 2 cm]
      \tikzset{VertexStyle/.style = {shape= circle, ball color= blue!20!white, inner sep= 2pt, outer sep= 0pt, minimum size = 2 pt}}
      \tikzset{EdgeStyle/.style   = {double= blue, double distance = 1pt}}

      \draw (0:4 cm) node[VertexStyle, font=\scriptsize] (0) {} ;
      \draw (45:4 cm) node[VertexStyle, font=\scriptsize] (1) {} ;
      \draw (90:4 cm) node[VertexStyle, font=\scriptsize] (2) {} ;
      \draw (135:4 cm) node[VertexStyle, font=\scriptsize] (3) {} ;
      \draw (180:4 cm) node[VertexStyle, font=\scriptsize] (4) {} ;
      \draw (225:4 cm) node[VertexStyle, font=\scriptsize] (5) {} ;
      \draw (270:4 cm) node[VertexStyle, font=\scriptsize] (6) {} ;
      \draw (315:4 cm) node[VertexStyle, font=\scriptsize] (7) {} ;
      
      \draw (10:2 cm) node[VertexStyle, font=\scriptsize] (s) {} ;
      \draw (55:2 cm) node[VertexStyle, font=\scriptsize] (s1) {} ;
      \draw (100:2 cm) node[VertexStyle, font=\scriptsize] (s2) {} ;
      \draw (145:2 cm) node[VertexStyle, font=\scriptsize] (s3) {} ;
      \draw (190:2 cm) node[VertexStyle, font=\scriptsize] (s4) {} ;
      \draw (235:2 cm) node[VertexStyle, font=\scriptsize] (s5) {} ;
      \draw (280:2 cm) node[VertexStyle, font=\scriptsize] (s6) {} ;
      \draw (325:2 cm) node[VertexStyle, font=\scriptsize] (s7) {} ;

      \draw [EdgeStyle] (0)--(s1) ;
      \draw [EdgeStyle] (0)--(s2) ;
      \draw [EdgeStyle] (0)--(s4) ;
      \draw [EdgeStyle] (1)--(s2) ;
      \draw [EdgeStyle] (1)--(s3) ;
      \draw [EdgeStyle] (1)--(s5) ;
      \draw [EdgeStyle] (2)--(s3) ;
      \draw [EdgeStyle] (2)--(s4) ;
      \draw [EdgeStyle] (2)--(s6) ; 
      \draw [EdgeStyle] (3)--(s4) ;
      \draw [EdgeStyle] (3)--(s5) ;
      \draw [EdgeStyle] (3)--(s7) ;
      \draw [EdgeStyle] (4)--(s5) ;
      \draw [EdgeStyle] (4)--(s6) ;
      \draw [EdgeStyle] (4)--(s) ;
      \draw [EdgeStyle] (5)--(s6) ;
      \draw [EdgeStyle] (5)--(s7) ;
      \draw [EdgeStyle] (5)--(s1) ;
      \draw [EdgeStyle] (6)--(s7) ;
      \draw [EdgeStyle] (6)--(s) ;
      \draw [EdgeStyle] (6)--(s2) ;
      \draw [EdgeStyle] (7)--(s) ;
      \draw [EdgeStyle] (7)--(s1) ;
      \draw [EdgeStyle] (7)--(s3) ;
    \end{tikzpicture}
  \end{center}
  \caption{Cay($D_8, \lsb sr, sr^2, sr^4 \rsb$)} \label{shuriken}
\end{figure}

\begin{ex}The graph in Figure~\ref{shuriken} is a regular graph with spectrum
\[ \lsb 3^{(1)}, \sqrt{3}^{(4)}, 1^{(3)}, \lp -1 \rp^{(3)}, \lp - \sqrt{3} \rp^{(4)}, \lp -3 \rp^{(1)} \rsb, \]
and the halved graphs have the spectrum
\[ \lsb 6^{(1)}, 0^{(4)}, (-2)^{(3)} \rsb. \]
Thus the halved graphs are strongly-regular and have the desired relationship, even though the graph itself has diameter four and six distinct eigenvalues.\end{ex}

In both these examples, \( G \) is regular. This turns out to be a necessary condition on such examples.

\begin{thm}\label{semiregularHalved}Let \( G = \lp B, C \rp \) be a \( \lp k, \ell \rp \) semiregular bipartite graph with \( k < \ell. \) Let \( H_B \) and \( H_C \) be the halved graphs with eigenvalues \( \theta_0 > \ldots > \theta_d \) and \( \tau_0 > \ldots > \tau_{d'}, \) respectively. If \( H_1, H_2 \) are distance-regular and there exist constants \( r,s \) such that the eigenvalues of \( G \) are
  \[ \bigcup_{i=0}^d \lsb \sqrt{r\theta_i+k} \rsb = \bigcup_{i=0}^{d'} \lsb \sqrt{s \tau_i+\ell} \rsb \bigcup \lsb 0^{(\abs{B}-\abs{C})} \rsb, \]
  then \( G \) is distance-biregular.
\end{thm}

\proof Note first that \( H_C \) cannot have more distinct eigenvalues that \( H_B, \) and therefore \( d' \leq d. \) Since \( H_B \) is distance-regular, we know that it has diameter \( d, \) and therefore \( G \) has diameter at least \( 2d. \) Then \( G \) has at least \( 2d+1 \) distinct eigenvalues, and since the eigenvalues of \( G \) are determined by the eigenvalues of \( H_B, \) we know \( G \) has at most \( 2d+2 \) distinct eigenvalues. Further, since \( k < \ell, \) we have \( \abs{B} > \abs{C}, \) and therefore 0 is an eigenvalue of \( G \). Since \( G \) is bipartite, it must therefore have an odd number of distinct eigenvalues, and thus \( G \) has \( 2d+1 \) distinct eigenvalues and diameter \( 2d \).

Since \( H_B \) is distance-regular, there exists a polynomial \( h^B \) of degree \( d \) such that
\[ h^{B} \lp A \lp H_B \rp \rp = A_d \lp H_B \rp. \]
Then
\[ p^{B} \lp x \rp := h^{B} \lp \frac{1}{r} x^2-k \rp \]
is polynomial of degree \( 2d \) satisfying
\[ p^{B} \lp A \lp G \rp \rp \verp{B}  = A_d \lp G \rp \verp{B} . \]

If \( d' = d, \) then we can define \( p^{C} \) analogously. Otherwise we have \( d' = d-1, \) and there is a polynomial \( p^{C} \) of degree \( 2d \) that is zero on the support of every eigenvalue in \( C \), and so
\[ p^{C} \lp A \lp G \rp \rp \verp{C} = \mathbf{0} = A_d \lp G \rp \verp{C}. \]
In either case, by Theorem~\ref{pBiDistance} we conclude that \( G \) is distance-regular.\qed

\section{Local Adjacency Polynomials}

For a vertex \( \vtxa \), the \( \vtxa \)-local $k$-adjacency polynomial \( Q_k \) is the unique polynomial of degree $k$ with \( \norm{Q_k}_{\vtxa} = 1 \) satisfying
\[ Q_k \lp \lambda \rp = \sup \lsb p \lp \lambda \rp : p \in \mathbb{R}_k \lb x \rb, \norm{p}_{\vtxa} \leq 1 \rsb. \]
The notion of local adjacency polynomials were introduced by Fiol and Garriga~\cite{garriga}, who showed that the \( \vtxa \)-local $k$-adjacency polynomial is indeed well-defined.

Let \( N_k \lp \vtxa \rp \) denote the set of vertices at distance at most \( k \) from \( \vtxa \). We combine and reformulate several results from~\cite{garriga} to obtain the following result.

\begin{prop}\label{bigQLittleQ}For all \( 0 \leq k \leq e, \) the \( \vtxa \)-local $k$-adjacency polynomial satisfies
  \[ Q_k \lp \lambda \rp \leq \frac{1}{v_{\vtxa}} \sqrt{\sum_{\vtxb \in N_k \lp \vtxa \rp} v_{\vtxa}^2}. \]
  Equality holds for some \( k \) if and only if there exists a polynomial \( q_k \) of degree \( k \) with \( \norm{q_k}_{\vtxa} = Q_k \lp \lambda \rp \) such that
  \[ q_k \lp A \rp \ver{\vtxa} = \sum_{\vtxb \in N_k \lp \vtxa \rp} \frac{v_{\vtxb}}{v_{\vtxa}} \ver{\vtxb}. \]
\end{prop}

\proof Let \( Q_k \) be the \( \vtxa \)-local $k$-adjacency polynomial, and consider
\[ \ip{Q_k \lp A \rp \ver{\vtxa}}{\mathbf{v}}. \]
By the spectral decomposition, we see this equals \( Q_k \lp \lambda \rp v_{\vtxa}. \)

On the other hand, we can view \( Q_k \lp A \rp \ver{\vtxa} \) as a vector which is indexed by vertices. Since \( Q_k \) is a polynomial of degree \( k, \) we know the entries of the vector indexed by vertices at distance at least \( k+1 \) from \( \vtxa \) must be zero. Let \( \mathbf{v}^k \) be the vector which is \( \mathbf{v}_{\vtxb} \) for the vertices \( \vtxb \in N_k \lp \vtxa \rp \) and 0 on the vertices at distance at least \( k+1 \) from \( \vtxa \). Then using Cauchy-Schwarz, we see
\[ \ip{Q_k \lp A \rp \ver{\vtxa}}{\mathbf{v}} \leq \norm{Q_k}_{\vtxa} \norm{\mathbf{v}^k} = \sqrt{\sum_{\vtxb \in N_k \lp \vtxa \rp} v_{\vtxb}^2}. \]

If equality holds, then we see that \( Q_k \lp A \rp \ver{\vtxa} \) must be a scalar multiple of \( \mathbf{v}^k \). In particular, to have the right norm, we must have
\[ Q_k \lp A \rp \ver{\vtxa} = \frac{1}{\sqrt{\sum_{\vtxc \in N_k \lp \vtxa \rp}}} \sum_{\vtxb \in N_k \lp \vtxa \rp} v_{\vtxb} \ver{\vtxb}. \]
Then if we let
\[ q_k \lp x \rp := Q_k \lp \lambda \rp Q_k \lp x \rp, \]
we can see it has the desired properties.

Conversely, suppose there exists a polynomial \( q_k \) of degree \( k \) with \( \norm{q_k}_{\vtxa} = Q_k \lp \lambda \rp \) such that
\[ q_k \lp A \rp \ver{\vtxa} = \sum_{\vtxb \in N_k \lp \vtxa \rp} \frac{v_{\vtxb}}{v_{\vtxa}} \ver{\vtxb}. \]
Then we can take the norm to see
\[ Q_k \lp \lambda \rp = \norm{q_k}_{\vtxa} = \frac{1}{v_\vtxa} \sqrt{\sum_{\vtxb \in N_k \lp \vtxa \rp} v_{\vtxb}^2}, \]
and so the earlier bound is tight.\qed

The polynomials \( q_0, \ldots, q_e \) can be defined in terms of the sequence of polynomials orthogonal to the local inner product, as was done by Fiol and Garriga~\cite{garriga}. This suggests a relationship between when the \( \vtxa \)-local adjacency polynomials achieve their maximum and when \( \vtxa \) is pseudo-distance-regular, which they were were able to prove. We reprove it here, using a different set-up.

\begin{thm}\label{excessPseudo}Let \( G \) be a connected graph with a vertex \( \vtxa \) of eccentricity \( e \) and \( \abs{\Phi_{\vtxa}} = e+1. \) Then \( G \) is pseudo-distance-regular around \( \vtxa \) if and only if
\[ Q_{e-1} \lp \lambda \rp =\frac{1}{v_{\vtxa}} \sqrt{\sum_{\vtxb \in N_{e-1}} v_{\vtxb}^2}. \]
\end{thm}

\proof If \( \vtxa \) is pseudo-distance-regular, then as discussed in Section~\ref{pseudo}, we know that there exist orthogonal polynomials \( p_0, \ldots, p_e \) such that \( p_i \lp A \rp \ver{\vtxa} \) is some scalar multiple of \( \mathbf{v} \) on the vertices at distance \( i \) from \( \vtxa \) and 0 otherwise. Thus we can take a linear combination of \( p_0, \ldots, p_{e-1} \) to get a polynomial \( q_{e-1} \lp x \rp \) of degree \( e-1 \) such that
\[ q_{e-1} \lp A \rp \ver{\vtxa} = \sum_{\vtxb \in N_{e-1} \lp \vtxa \rp} \frac{v_{\vtxb}}{v_{\vtxa}} \ver{\vtxb}. \]
Letting
\[ Q_{e-1} \lp x \rp := \frac{q_e \lp x \rp}{\norm{q_e}_u}, \]
we see this is a polynomial of degree \( e-1 \) and norm 1 satisfying
\[ Q_{e-1} \lp \lambda \rp =\frac{1}{v_{\vtxa}} \sqrt{\sum_{\vtxb \in N_{e-1}} v_{\vtxb}^2}, \]
so it is indeed the \( \lp e-1 \rp \)-excess polynomial.

Conversely, suppose
\[ Q_{e-1} \lp \lambda \rp =\frac{1}{v_{\vtxa}} \sqrt{\sum_{\vtxb \in N_{e-1}} v_{\vtxb}^2}. \]
By Proposition~\ref{bigQLittleQ}, we see that there is a polynomial \( q_{e-1} \) of degree \( e-1 \) such that
\[ q_{e-1} \lp A \rp \ver{\vtxa} = \sum_{\vtxb \in N_{e-1} \lp \vtxa \rp} \frac{v_{\vtxb}}{v_{\vtxa}} \ver{\vtxb}. \]

Let \( \theta_0, \ldots, \theta_e \) be the eigenvalues in the support of \( \vtxa \). Define
\[ q_e \lp x \rp = \frac{1}{v_{\vtxa}^2} \prod_{\theta_r \neq \lambda} \frac{x-\theta_r}{\lambda-\theta_r}. \]
Note that \( q_e \) is polynomial of degree \( e \) such that
\[ q_e \lp A \rp \ver{\vtxa} = \frac{1}{v_{\vtxa}} \mathbf{v} = \sum_{\vtxb} \frac{v_{\vtxb}}{v_{\vtxa}} \ver{\vtxb}. \]

Then \( q_e - q_{e-1} \) is a polynomial of degree at most \( e \) such that \( \lp q_e- q_{e-1} \rp \lp A \rp \ver{\vtxa} \) is nonzero precisely on the vertices of \( \vtxa \) at distance \( e \) from \( \vtxa \). By Theorem~\ref{pseudoLast}, we see that \( \vtxa \) is pseudo-distance-regular.\qed

This alternate characterization of pseudo-distance-regular vertices can be applied globally to give us a characterization of when a graph is distance-regular or distance-biregular. The next result is due to Fiol and Garriga~\cite{garriga} though the formalization is different. The proof here gives a preview of the main ideas before we get into the casework of the distance-biregular analogue.

\begin{thm}[Fiol and Garriga~\cite{garriga}]\label{pDistance}Let \( G \) be a regular connected graph with vertex set \( V \) diameter \( d \), and largest eigenvalue \( \lambda \). Then \( G \) is distance-regular if and only if it has \( d+ 1 \) distinct eigenvalues and there exists a polynomial \( p\) of degree \( d \) such that
  \[ \norm{p}_G^2 =  p \lp \lambda \rp = \abs{\lsb \vtxa \in V : d \lp \vtxb, \vtxa \rp = d \rsb}\]
  for all vertices \( \vtxa. \)
\end{thm}

\proof If \( G \) is distance-regular, then it has \( d+1 \) distinct eigenvalues and the distance polynomial \( p_d \) is the desired polynomial.

Let \( n \) be the number of vertices of \( G, \) and let \( t \) be the number of vertices at distance \( d \) from any vertex in the graph. Define
\[ q \lp x \rp = n \prod_{\theta_r \neq \lambda} \frac{x- \theta_r}{\lambda - \theta_r}. \]
Note that \( \norm{q}_G^2 = n \) and \( \ip{q}{p}_G = t \). Thus, the projection of \( q \) onto \( p \) is
\[ \frac{\ip{q}{p}_v}{\norm{p}_G^2} p = p. \]

We can extend \( p \) into an orthogonal sequence of polynomials \( p_0, \ldots, p_d = p \). Since \( q-p \) is orthogonal to \( p, \) we must be able to write it as a linear combination of \( p_0, \ldots, p_{d-1} \), therefore
\[ Q \lp x \rp := \frac{1}{\sqrt{n-t}} \lp q \lp x \rp - p \lp x \rp \rp\]
is a polynomial of degree at most \( d-1 \). Note that
\[ \norm{Q}_G^2 = \frac{1}{n-t} \lp \norm{q}_G^2 - 2 \ip{q}{p}_G + \norm{p}_G^2 \rp = 1. \]

We have
\[ Q \lp \lambda \rp = \frac{1}{\sqrt{n-t}} \lp n-t \rp = \sqrt{n-t} = \sqrt{n} \sqrt{\lp n-t \rp \frac{1}{n}}. \]
Since \( G \) is regular, the Perron eigenvector is \( \frac{1}{\sqrt{n}} \mathbf{1}, \) and therefore for any vertex \( \vtxa, \) we must have
\[ Q \lp \lambda \rp = \frac{1}{v_{\vtxa}} \sqrt{\sum_{\vtxb \in N_{d-1} \lp \vtxa \rp} v_{\vtxb}^2}. \]
Thus \( \norm{Q}_{\vtxa} \geq 1, \) or else \( Q \) would exceed the bound in Proposition~\ref{bigQLittleQ}.

It follows
\[ n \leq \sum_{\vtxa \in V} \norm{Q}_{\vtxa}^2 = \sum_{\vtxa \in V} \sum_{r=0}^d Q^2 \lp \theta_r \rp \lp E_r \rp_{\vtxa, \vtxa} = n \sum_{r=0}^d Q^2 \lp \theta_r \rp \frac{1}{n} \mathrm{tr} \lp E_r \rp = n \norm{Q}_G^2 =1. \]
Equality holds, therefore \( \norm{Q}_{\vtxa} =1 \) for all vertices \( \vtxa \). Then since \( Q \) is a polynomial of \( \vtxa \)-local norm 1 attaining the bound  of Theorem~\ref{excessPseudo}, we have that \( \vtxa \) is pseudo-distance-regular. This is true of every vertex in \( G, \) and \( G \) is regular, therefore \( G \) is distance-regular.\qed

\begin{thm}\label{pBiDistance}Let \( G = \lp B, C \rp \) be a connected semiregular bipartite graph with diameter \( d \) and \( d+1 \) distinct eigenvalues. Then \( G \) is distance-biregular if and only if there exist polynomials \( p^{B} \) and \( p^{C} \) of degree \( d \) such that, for every vertex \( \vtxb \in B \), we have
\[ \norm{p^B}_B^2 =  p^{B} \lp \lambda \rp = \abs{\lsb \vtxa \in V : d \lp \vtxb, \vtxa \rp = d \rsb}\]
and, for every vertex \( \vtxc \in C, \) we have
\[ \norm{p^C}_C^2 =  p^{C} \lp \lambda \rp = \abs{\lsb \vtxa \in V : d \lp \vtxc, \vtxa \rp = d \rsb}. \]
\end{thm}

\proof This is clearly true if \( G \) is distance-biregular, so we may suppose there exist polynomials \( p^{B}, p^{C} \) of degree \( d \) such that
\[ \norm{p^B}_B^2 =  p^{B} \lp \lambda \rp = \abs{\lsb \vtxa \in V : d \lp \vtxb, \vtxa \rp = d \rsb} =: t\]
and, for every vertex \( \vtxc \in C, \) we have
\[ \norm{p^C}_C^2 =  p^{C} \lp \lambda \rp = \abs{\lsb \vtxa \in V : d \lp \vtxc, \vtxa \rp = d \rsb}. \]

Assume without loss of generality that vertices in \( B \) have maximum eccentricity \( d \), and let \( d' \) be the maximum eccentricity of vertices in \( C \). If \( d \) is odd, then by Lemma~\ref{biRegularOddDiameter}, \( G \) is regular, and Theorem~\ref{pDistance} says \( G \) is a bipartite distance-regular graph. Thus we may assume \( d \) is even. In particular, the vertices at distance \( d \) from vertices in \( B \) all lie in \( B \).

Let
\[ q^{B} \lp x \rp = 2 \abs{B} \prod_{\theta_r \neq \lambda} \frac{x-\theta_r}{\lambda-\theta_r}. \]
We have \( \ip{q^{B}}{p^{B}}_B = t, \) so the projection of \( q^{B} \) onto \( p^{B} \) is \( p^{B} \). In particular, extending \( p^{B} \) to an orthogonal sequence of polynomials we see
\[ Q^{B} \lp x \rp := \frac{1}{\sqrt{2 \abs{B}-t}} \lp q^{B} \lp x \rp - p^{B} \lp x \rp \rp \]
must have degree at most \( d-1 \).

Note
\[ \norm{Q^{B}}_B^2 = \frac{1}{2 \abs{B}-t} \lp q^{B} \lp \lambda \rp^2 \frac{1}{2 \abs{B}} -2  q^{B} \lp \lambda \rp p^{B} \lp \lambda \rp \frac{1}{2 \abs{B}} + \norm{p^{B}}_B^2 \rp = 1. \]

For \( \vtxb \in B, \) the bound of Proposition~\ref{bigQLittleQ} for the \( \vtxb \)-local $d$-excess polynomial is
\[ \sqrt{2 \abs{B}} \sqrt{\sum_{\vtxc} v_{\vtxc}^2} = \sqrt{2 \abs{B}} \sqrt{\frac{1}{2 \abs{B}} \lp\abs{B}-t+\abs{C} \frac{\ell}{k} \rp} = \sqrt{2 \abs{B}-t}. \]
Since
\[ Q^{B} \lp \lambda \rp = \frac{1}{\sqrt{2 \abs{B}-t}} \lp q^{B} \lp \lambda \rp - p^{B} \lp \lambda \rp \rp = \sqrt{2 \abs{B}-t}, \]
we therefore have that
\[ \norm{Q^B}_{\vtxb} \geq 1. \]
Then
\[ \abs{B} \leq \sum_{\vtxb \in B} \norm{Q^{B}}_{\vtxb} = \sum_{r=0}^d Q^{B} \lp \theta_r \rp^2 \sum_{\vtxb \in B} \lp E_r \rp_{\vtxb} = \abs{B} \norm{Q^{B}}_B =\abs{B}. \]
Since equality holds, we have \( \norm{Q^{B}}_{\vtxb} = 1 \) for all \( \vtxb \in B, \) and so by Theorem~\ref{excessPseudo}, we have that every vertex in \( \vtxb \) is pseudo-distance-regular. In fact, \( p^{B} \) is the \( d \)-local polynomial for all \( \vtxb \in B, \) so the set \( B \) is pseudo-distance-regular.

If \( d' = d, \) then the same argument holds for \( C. \) Otherwise, we know \( p^{C} \) is a polynomial of degree \( d \) that is zero on the eigenvalues in the support of \( C, \) so every vertex \( \vtxc \in C \) has eccentricity \( d-1 \) and \( d \) eigenvalues in \( \Phi_{\vtxc}. \) Since the set \( B \) is pseudo-distance-regular, we know there exists a polynomial \( p_{d-1}^{B} \) of degree \( d-1 \) such that \( p_{d-1}^{B} \lp A \rp \verp{B} \) is nonzero precisely on the vertices at distance \( d-1 \) from \( B \). Since \( p_{d-1}^{B} \lp A \rp \) is symmetric and \( d-1 \) is odd, it follows that \( p_{d-1}^{B} \lp A \rp \verp{C} \) is nonzero precisely on the vertices at distance \( d-1 \) from vertices in \( C \). By Theorem~\ref{pseudoLast}, the set \( C \) is also pseudo-distance-regular. Since every vertex in \( G \) is pseudo-distance-regular and \( G \) is bipartite, by Theorem~\ref{pseudoGlobal} we know that \( G \) is distance-biregular.\qed

\section{Spectrum of a Distance-Biregular Graph}

One way that the spectral excess characterization of Fiol and Garriga~\cite{garriga} has been used has been in studying graphs with the spectrum of a distance-regular graph. Just because a graph \( G' \) has the same spectrum as a distance-regular graph does not mean that \( G' \) is distance-regular. Hoffman~\cite{hoffmanPolynomial} gave the example of a graph cospectral to the Hamming graph \( H \lp 4, 2 \rp \) which was not distance-regular. Since \( H \lp 4, 2 \rp \) is bipartite, this also shows that a graph cospectral to a distance-biregular graph need not be distance-biregular.

It is worth noting that a graph cospectral to a connected regular graph will be regular, with the valency equal to the largest eigenvalue, so the valency of the graph is spectral information. The same is not true of semiregular bipartite graphs. The graphs \( K_{1,4} \) and \( C_4 \cup K_1 \) are cospectral, even though they are not both connected and they have different valencies. Thus when we are considering semiregular bipartite graphs with a particular spectrum, we will allow ourselves the additional spectral information of the valencies of the vertices in each set of the partition.

Haemers~\cite{haemersSpectrumDistance} gave more examples of non-distance-regular graphs cospectral to distance-regular graphs, as well as studying what conditions were enough to guarantee a graph was distance-regular from the spectrum. Subsequently, Van Dam and Haemers~\cite{vanSpectral} provided more conditions to prove a graph with spectrum of a distance-regular graph was distance-regular. Abiad, Van Dam, and Fiol~\cite{abiad} further improved characterizations of when the spectrum forced a graph to be distance-regular.

Some of these characterizations can be derived using Theorem~\ref{pDistance} of Fiol and Garriga~\cite{garriga}. For instance, Van Dam and Haemers~\cite{whichGraphs} used Theorem~\ref{pDistance} to prove that a graph cospectral to a distance-regular graph of diameter \( d \) and girth at least \( 2d-1 \) is distance-regular. This result was first shown by Brouwer and Haemers~\cite{gewirtz}, and if the graph is bipartite, Van Dam and Haemers~\cite{vanSpectral} improved this to prove that a graph cospectral to a bipartite distance-regular graph of diameter \( d \) and girth at least \( 2d-2 \) is distance-regular. The equivalent result for distance-biregular graphs can be derived by Theorem~\ref{pBiDistance}.

\begin{thm}\label{cospectralGirth}Let \( G \) be a distance-biregular graph of diameter \( d \) and girth at least \( 2d-2 \) and valencies \( k, \ell. \) Let \( G' \) be a \( \lp k, \ell \rp \)-semiregular graph cospectral to \( G. \) Then \( G' \) is distance-biregular.
\end{thm}

\proof The girth is determined by the spectrum, so \( G' \) has girth \( 2d-2 \), and since \( G' \) only has \( d+1 \) distinct eigenvalues, \( G' \) can only have diameter \( d-1 \) or \( d \). If \( G' \) has diameter \( d-1, \) it is a generalized polygon, which Godsil and Shawe-Taylor~\cite{distanceRegularised} proved was distance-biregular. Otherwise, we know that \( G' \) has diameter \( d \) and girth \( 2d-2. \)

Let \( B \) be the set of vertices of valency \( k \) in \( G' \) and \( C \) be the set of vertices of valency \( \ell. \) Suppose without loss of generality that the vertices in \( B \) have eccentricity \( d \), and fix \( \vtxa \in B \). We wish to count the number of vertices at distance \( d \) from \( \vtxa \).

Suppose \( d \) is odd, and let \( 0 < i < d \) be another odd integer. Then since \( g \geq 2d-2, \) there can only be one walk from \( \vtxa \) to a vertex at distance \( i \) from \( \vtxa. \) Therefore, the number of vertices at distance \( i \) from \( \vtxa \) must be equal to the number of paths of non-backtracking walks of length \( i \), which was shown in~\cite{solo} to be
\[ k \lp \ell-1 \rp^{\frac{i-1}{2}} \lp k-1 \rp^{\frac{i-1}{2}}. \]
Thus, the number of vertices at distance \( d \) from \( \vtxa \) is
\[ k_d^{B} := \abs{C}-\sum_{i=0}^{\frac{d-1}{2}} k \lp \ell-1 \rp^{i} \lp k-1 \rp^{i}. \]
Similarly, if \( d \) is even and \( 1 < i < d \) is another even integer, the number of vertices at distance \( i  \) from \( \vtxa \) must be
\[ k_d^{B} := \abs{B} - \sum_{i=1}^{\frac{d}{2}-1} k \lp \ell-1 \rp^{i} \lp k-1 \rp^{i-1} - 1. \]
In both these cases, \( k_d^{B} \) is the same as the number of vertices at distance \( d \) from \( \vtxa \) in \( G, \) so Theorem~\ref{pBiDistance} tells us that there is a polynomial \( p^{B} \) of degree \( d \) such that
\[ \norm{p^B}_B^2 =  p^{B} \lp \lambda \rp = k_d^{B} \]
where \( \lambda \) is the largest eigenvalue. A similar argument holds for a vertex in \( C, \) so we have \( k_d^{C} \) is the number of vertices at distance \( d \) from any given vertex in \( C, \) and a polynomial \( p^{C} \) of degree \( d \) such that
\[ \norm{p^{C}}_C^2 = p^{C} \lp \lambda \rp = k_d^{C}, \]
and therefore by Theorem~\ref{pBiDistance} we know \( G' \) is distance-biregular.\qed

Partial geometries are distance-biregular graphs with girth six and diameter four. Therefore, any graph with the spectrum of a partial geometry is distance-biregular. This extends the result of Van Dam and Haemers~\cite{vanSpectral} that any graph the spectrum of a regular partial geometry is distance-regular.

In fact, a stronger version of Theorem~\ref{cospectralGirth} is possible. Abiad, Van Dam, and Fiol~\cite{abiad} proved that any regular bipartite graph with \( d+1 \) distinct eigenvalues and girth at least \( 2d-2 \) is distance-regular. The author~\cite{solo} proved an analogous result for semiregular bipartite graphs. However, for the question of determining when a graph with the spectrum of a distance-biregular graph is distance-biregular, Theorem~\ref{cospectralGirth} is sufficient, and the proof demonstrates how using Theorem~\ref{pBiDistance} can be used to prove distance-biregularity from the spectrum of certain graphs.

\section*{Acknowledgements}
I would like to thank Chris Godsil for his help at every stage of the process.

%\bibliographystyle{acm}
%\bibliography{../ExtremeGhostbusting/Bibliography}

\end{document}